# Optimization of waste collection through the sequencing of micro-routes and transfer station convenience analysis: an Argentinian case study

**SOFÍA A. MOLFESE GRECO[1] - DIEGO G. ROSSIT[1][2]* - MARIANO FRUTOS[1][3] - ANTONELLA CAVALLIN[1]**

[1] Departamento de Ingeniería, Universidad Nacional del Sur

[2] Instituto de Matemática de Bahía Blanca UNS-CONICET

[3] Instituto de Investigaciones Económicas y Sociales del Sur UNS-CONICET

*Corresponding autor: diego.rossit@uns.edu.ar – 1253 Alem Avenue, Bahía Blanca B8000CPB, Argentina*

## ABSTRACT

Municipal solid waste management is a paramount activity in modern cities due to the environmental, social and economic problems that can arise when mishandled. In this work, the sequencing of micro-routes in the Argentine city of Bahía Blanca is addressed, which is modeled as a vehicle routing problem with travel time limit and the vehicle's capacity. Particularly, we propose two mathematical formulations based on mixed-integer programming and we solve a set of instances of the city of Bahía Blanca based on real data. Moreover, with this model we estimate the total distance and travel time of the waste collection and use this data to analyze the possibility of installing a transfer station. The results demonstrate the competitiveness of the approach to resolve realistic instances of the target problem and suggest the convenience of installing a transfer station in the city considering the reduction of the traveled distance.

## 1. INTRODUCTION

In modern cities, waste management has become a critical issue for enhancing the sustainability and improving the quality of life of citizens (De Souza et al., 2017).

Moreover, this issue becomes even more critical if it is considered that the amount of waste generation rate per capita will continue increasing in the following decades (Hoornweg and Bhada, 2012; Hoornweg et al., 2015).

The reverse logistic chain of Municipal Solid Waste (MSW) is a special case of a reverse supply chain that is defined as "a network consisting of all entities involved in the flow of disposed products leaving the point of consumption. It includes collection, transportation, recovery and disposal of waste. Its purpose is to recapture or create value and/or proper disposal" (Van Engeland et al., 2020). This subject has received a large interest from both professionals and academic communities in the last decades. An efficient provision of MSW management services is a key element to diminish the environmental impact of human activities on the environment, through a responsible treatment and/or final disposition of the generated waste, and to postpone the depletion of limited resources, through the recovery of reusable resources from waste (Das et al., 2019; Rossit and Nesmachnow, 2022).

The problem that can be caused by inadequate management of the large volumes of the waste generated in various urban agglomerations is well known (Poletto et al., 2016). This has not only environmental implications, but also economic and social ones. Developing countries, such as Argentina, that usually lack large budgetary resources and/or advanced decision support tools, generally have limited tools to reduce the negative impact of the MSW system on the quality of life of citizens (Abarca et al., 2013). Among the diverse stages in the reverse logistic chain of municipal waste the collection phase is considered one of the most costly (Das and Bhattacharyya, 2015). Thus, providing computer aid decision support tools to improve the decision-making process in this collection stage can help to enhance efficiency of the overall MSW system and reduce the associated costs. This is specially important for Argentina, a country which has a remarkably high logistic and transport costs (Musante, 2021). In the Argentine city of Bahía Blanca more than 157 thousand tons of MSW are collected per year[1]. Moreover, the MSW system, including the collection and the landfill management, consumes approximately 20% of the budget of the municipal government[2]. Currently, in the city waste collection is carried out with a door-to-door system in which the vehicle that collects the waste must visit all the dwellings of the city every day of the week except on Sunday. Collection vehicles carry the waste from the dwelling directly to the landfill which is located outside the city. The travel distance of the round trip between the downtown of the city and the landfill is about 34 kilometers with an estimated travel time of 45 minutes. This situation is common in MSW systems since environmental regulations and citizens' concern have gradually forced landfills to be located further and further from the downtown (Ghiani et al., 2021).

## 1.1. Contribution and organization of this work

In this work we aim to make two main contributions. Firstly, we propose two mixed integer programming models based on the classical vehicle routing problem to model two different situations: the current situation of the city with no transfer station in the MSW system and the hypothetical situation in which a transfer station is included in the system. Both models aim at minimizing the total traveled distance subjected to different operational constraints such as capacity of the vehicle and maximum travel time. Although these models are tailor-made for the city of Bahía Blanca, they can also be applied to other cities with similar characteristics. Secondly, we analyze the convenience of the inclusion of the transfer station -studying two different locations that were proposed in accordance with practitioners- and the impact on the MSW in terms of total traveled distance. Transfer stations usually work as intermediate facilities in which the collection vehicles deposit the MSW gathered from the city. They are generally located nearer to the urban areas than the landfill or final disposition plants for temporarily holding the MSW (Yadav and Karmakar, 2020). Then, when enough waste is accumulated, the MSW is transferred by larger or transfer vehicles, which usually have a smaller operating cost than the common collection vehicles (Höke and Yalcinkaya, 2021; Eshet et al., 2007; Kirca and Erkip, 1988; Rathore and Sarmah, 2019), to the final destination. Depending on the level of technification, transfer stations can apply some processing to the MSW. For example, further compacting or classification in different streams for recycling. According to Kirca and Erkip (1988), depending on the particular characteristics of the city, transfer stations can have the following benefits for the overall MSW system: i) the collection vehicles spend less time traveling and more time in collecting waste, ii) the personnel is more effectively employed since more time is spent on collection, iii) the response rate to unexpected service emergencies from the collection vehicles is reduced since vehicles are usually near to the area in which they operate, iv) depending on the level of technification, transfer stations can apply some processing to the MSW such as waste sorting and classification which can reduce the material that is finally sent to the final destination, v) and roads to landfill or processing plants are less congested as one transfer vehicle may replace several collection vehicles (since increasing the size of the vehicle creates less congestion than increasing the number of vehicles). Additionally, Rathore and Sarmah (2019) state that transfer stations can contribute to a reduction in pollution impact of the system and the quality of service provided to the citizens.

---

[1] https://gobiernoabierto.bahia.gob.ar/ambientales/residuos-domiciliarios/
[2] https://www.bahia.gob.ar/economia/presupuesto/

This work is organized as follows. In the rest of this Section the main related works are presented. In Section 2 the case study and the research methodology are described. In Section 3 the main results of the computational experimentation are presented and analyzed. Finally, Section 4 outlines the main conclusions along with the future research lines.

## 1.2 Related works

Regarding transfer station analysis, Höke and Yalcinkaya (2021) presented a geographic information system (GIS)-based analysis to study the suitability of different potential locations of transfer stations and waste collection routes for various collection vehicle capacities in the Turkish city of İzmir. The simulated scenarios with addition of the transfer station allowed to reduce the collection time by 9%. Rathore and Sarmah (2019) presented an analysis for locating a transfer station in a case study in Bilaspur, India with a GIS-based tool. They found that the inclusion of a transfer station in the city allowed to reduce the cost of the MSW system by 30%. Then, a set of works studied the potential location of transfer stations using transport models instead of vehicle routing problems, i.e., considering direct convoys from municipalities to transfer stations (which is a simpler optimization model), in larger regions such as Hellenic region, Greece (Chatzouridis and Komilis, 2012) and New South Wales, Australia (Asefi et al., 2015, 2019). Jia et al. (2022) analyzed the location of the transfer plant within a transport problem in which they considered the adjustment of the previous existing facility in a case study of Beijing, China. Vargas et al. (2022) performed a study in the Colombian city of Medellin in which they found that the introduction of transfer stations allowed to reduce the size of the vehicle collection fleet.

In Argentina, the design of MSW collection routes is based mostly on empirical knowledge of decision makers and not in systematized procedures or decision aid tools (Cavallin et al., 2020). The reason of the relatively poor application of systematized procedures in MSW management in Argentinian municipalities -which are in charge of MSW- can be connected to several obstacles (Cavallin, 2019): i) poor documented data to make informed decisions, e.g., generation rate per neighborhood, recycling rates, waste composition and information about the costs of the system; ii) lack of budgetary resources and regulatory framework that incentivize the investment in technology in

MSW system; iii) not developed markets for recycled product which do not promote the development of more technified MSW system that allow recuperation of resources from waste. However, -although not directly connected to the convenience transfer station analysis- some contributions have been made in different Argentine cities related to improved MSW logistics. In Buenos Aires, Bonomo et al. (2012) worked on optimizing routes with the aim of not only minimizing travel distances but also of minimizing vehicle wear. In Morón, Braier et al. (2017) developed an integer programming model for recyclable waste collection. Regarding the particular case of Bahía Blanca, some studies have been performed to improve the efficiency of the MSW system. A set of works addressed the first stage of the reverse logistic chain of MSW focusing on the storage level for accumulating the waste generated at households. Among the diverse storage levels that can be implemented (Gallardo et al., 2015; Rossit and Nesmachnow, 2022), in Cavallin et al. (2020) a proposal to migrate from the current door-to-door system to a waste bins-based system was analyzed as a valid approach to improve the efficiency of the MSW of the system. Thus, several works addressed the location of waste bins among the city with multiobjective exact (Rossit et al., 2017) and heuristic methods (Toutouh et al., 2018, 2020). Moreover, some works explicitly consider the location of garbage accumulation points taking into account a proxy to posterior routing costs among the different objectives (Rossit et al., 2018, 2020). Other works addressed the waste collection problem considering the previous location of waste bins (Fermani et al., 2020; Rossit et al., 2021) while others considered both problems simultaneously, i.e., optimizing the sizing of waste bins and designing the collection routes (Mahéo et al., 2020, 2022). Finally, in Vazquez et al. (2020) an analysis of the composition and amount of waste in Bahía Blanca was performed based on thorough field work.

This work addressed MSW collection from the perspective of considering fixed micro-routes. This is not a common issue in the related literature but it has certain benefits such as not affecting to a large extent the work routine of drivers and practitioners -which can be reluctant to drastic changes- and reduces the complexity of the problem allowing to address the routing plan of a whole city in a unique model. Thus, we consider that in the related bibliography there is still room to propose contributions that address these aspects. Moreover, this work is based on a case study of Argentina, where MSW systems usually do not take advantage from decision support tools and are only based on empirical knowledge of decision makers. Without underestimating the invaluable contribution that can be made from practice and experience, the bibliography is extensive in the contribution both in terms of cost reduction and environmental impact that the use of decision support techniques can provide. In order to enhance MSW system efficiency, we consider that providing tools that help to optimize waste collection is promising since the country has one of the largest logistic costs in the region (Musante, 2021).

## 2. MATERIALS AND METHODS

This section presents the description of the study area and the methodology used for the analysis.

### 2.2 Definition and description of the study areas

The city currently has around 300,000 inhabitants and is an important university, commercial and industrial center in the South of Argentina, having one of the main export ports of the country. The average population density of the city is about 2300 hab/km$^2$. Currently, waste collection is performed using a door-to-door system in which the collection vehicle visits each dwelling -which have their own waste bin in the sidewalk close to the kerbside- to collect the (unclassified) waste every day (Cavallin et al., 2020).

Waste collection in the city is organized in different districts or zones, so-called micro-routes. A micro-route is a relatively compact area in which the collection is performed by the same vehicle without intermediate stops. Fig. 1a shows the 32 micro-routes that are currently used in Bahía Blanca. From this total, 17 micro-routes correspond to the day shift and 15 micro-routes to the night shift. In general, the detailed internal path of the vehicle inside the micro-route area is already established. In these systems based on micro-routes, decision-makers and personnel are usually reluctant to modify the micro-routes covered area or the internal path (Vargas et al., 2021). However, the internal route of a micro-route is quite optimized in terms of distance and travel times given the experience of the personnel in charge (see, e.g., the internal path of micro-route number 22 in Fig. 1b). Moreover, the zones are typically compact areas -as least as compact as the urban street network allows- which has benefits of the operative conditions of the routes as has been stated in Rossit et al. (2019). Thus, taking into account these reasons, in this work, the micro-routes are considered as fixed and the problem addressed to improve collection logistics is their sequencing, that is, which micro-routes and in which order are visited by each collection vehicle.

The waste that is collected in the city is directly taken to the landfill (Location L in Fig. 1a). Thus, part of this work is to study the convenience of installing a transfer plant. For these purposes, two locations are considered: Location A, which corresponds to the current depot and administrative location of the collection company, and Location B, which is an industrial area and was considered a suitable place for this analysis by practitioners.

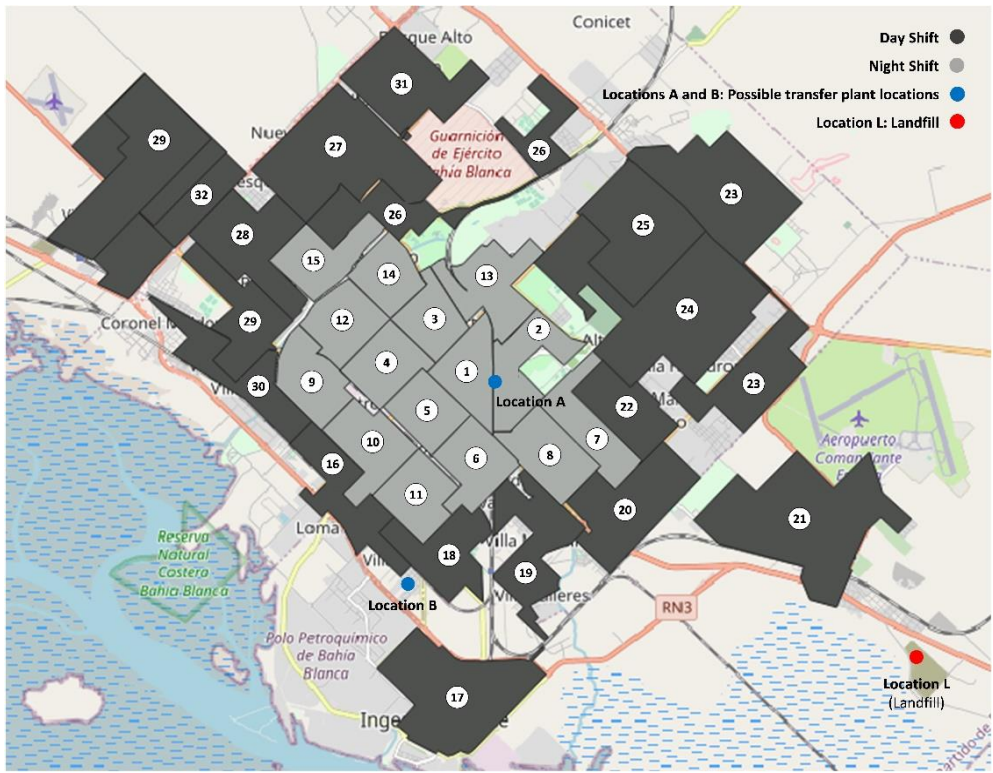

**a)** Micro-routes of the day and night shifts of Bahía Blanca.

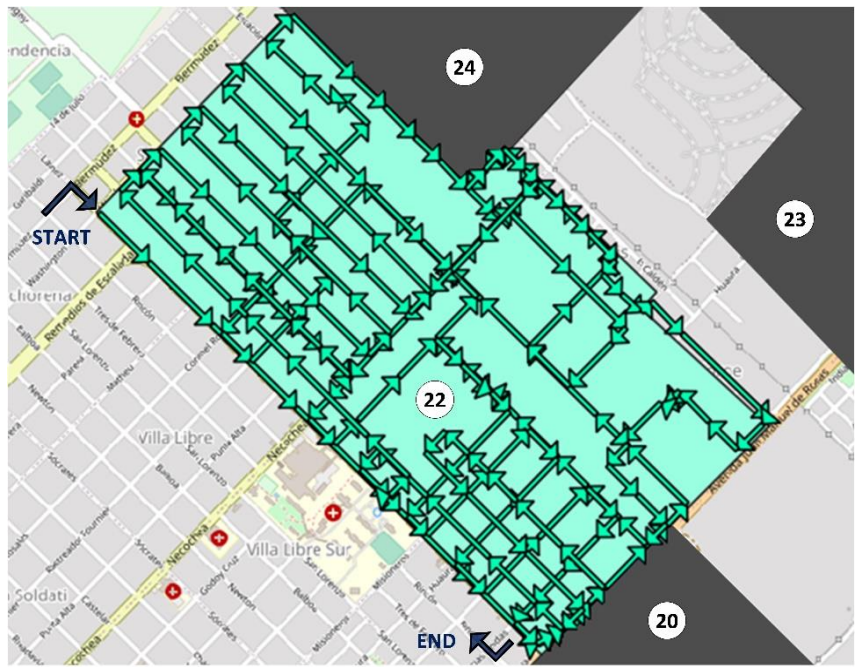

**b)** Internal route of the micro-route 22 of the day shift.

**Figure 1**. Micro-routes of Bahía Blanca.

## 2.3 Methodology

The problem of collecting the waste generated in each micro-route can be modeled as a vehicle routing problem with capacity and time limit constraints. We present two models. One for the current situation which involves that the collection vehicles carry the waste directly to the landfill, hereafter *current situation* (CS) case, and another one for the simulated study of including a transfer station, hereafter *transfer station* (TS) case. We present the two models below.

### 2.3.1 Mathematical formulation for the current situation case.

In the current situation of Bahía Blanca, which corresponds to the CS case, the depot of the fleet of vehicles is located in the site of the company providing the waste collection service. Therefore, the routes have to start and end in this location. Additionally, vehicles have to unload waste at the landfill. This occurs when the capacity of the vehicle is near to be full (and it cannot visit any more micro-routes before unloading) or it has to finish the tour due to the end of the shift -since for sanitary reasons it cannot return to the depot loaded with waste-. Taking these aspects into account, for the CS case, we introduced a model based on the three-index (vehicle flow) formulation for Capacitated Vehicle Routing Problem (CVRP) in which the set of vehicles (or routes) are identified (Irnich et al., 2014). Additionally, to the basic model CVRP we add the consideration of a time limit for the routes and the visit to the landfill as an intermediate facility (Hemmelmayr et al., 2013; Markov et al., 2015). In this sense, all the vehicles start and end at the depot empty. Before going back to the depot, vehicles have to visit the landfill in order to unload the waste. Moreover, the vehicles can visit the landfill more than once per route in order to unload and visit more micro-routes. Therefore, the model for the CS case can be formulated as follows. The following sets are considered: the set $M$ of micro-routes, the set $M' = M \cup \{0\} \cup \{l\}$ in which $0$ represents the depot and $l$ represents the landfill, and the set $K$ of vehicle routes. The following variables are defined: binary variable $x_{ijk}$ that takes the value 1 if vehicle $k \in K$ uses the path from the exit of micro-route $i \in M'$ to the entrance of micro-route $j \in M'$ and it takes the value 0 otherwise; and the positive continuous variable $v_{ijk}$ that indicates the load of vehicle $k \in K$ in the path from micro-route $i \in M'$ to micro-route $j \in M'$. Then, the parameter $Q$ is defined as the capacity of each collection vehicle, $d_i$ is the distance between the exit of micro-route $i \in M'$ and the entrance of micro-route $j \in M'$, $q_i$ is the quantity of waste to be collected on micro-route $i \in M$, $s_i$ is the collection time on micro-route $i \in M$, $h_i$ is the time it takes a vehicle to go from exit of micro-route $i \in M'$ and the entrance of micro-route $j \in M'$ and $T$ is the time limit for any route. Thus, the model is defined as a mixed integer programming model by Equations (1) - (11).

$$min \sum_{i,j \in M', \forall j \neq i, k \in K} x_{ijk} d_{ij} \qquad (1)$$

Subject to:

$$\sum_{j \in M', j \neq i, k \in K} x_{ijk} = 1, \forall i \in M \qquad (2)$$

$$\sum_{i \in M \cup \{l\}} x_{0ik} \leq 1, \forall k \in K \qquad (3)$$

$$\sum_{j \in M} x_{j0k} = 0, \forall k \in K \qquad (4)$$

$$\sum_{i \in M', j \neq i} x_{ijk} - \sum_{i \in M', j \neq i} x_{jik} = 0, \forall j \in M', k \in K \qquad (5)$$

$$v_{lik} = 0, \forall i \in M \cup \{0\}, k \in K \qquad (6)$$

$$v_{0ik} = 0, \forall i \in M \cup \{l\}, k \in K \qquad (7)$$

$$\sum_{i \in M \cup \{0\}, j \neq i} v_{ijk} + q_j \leq \sum_{i \in M \cup \{l\}, j \neq i} v_{jik} + Qt \left(1 - \sum_{i \in M \cup \{0\}, j \neq i} x_{ijk}\right), \forall j \in M, k \in K \qquad (8)$$

$$v_{ijk} \leq Qt \, x_{ijk}, \forall \, i, j \in M, j \neq i, k \in K \qquad (9)$$

$$x_{ijk} \leq \sum_{g \in M} x_{0gk}, \forall \, i, j \in M \cup \{l\}, j \neq i, k \in K \qquad (10)$$

$$\sum_{i, j \in M', j \neq i} (s_j + h_{ij}) \, x_{ijk} \leq T, \forall k \in K \qquad (11)$$

$$x \in B, v \geq 0$$

The proposed objective consists of minimizing the total distance traveled and is expressed in equation (1). Equation (2) sets that every micro-route has to be visited. Equation (3) establishes that a vehicle can only leave the depot at most once. Equation (4) forces that a vehicle visits the landfill to dump the waste before going back to the depot. Equation (5) sets that if a micro-route is visited by vehicle, the vehicle must also leave the micro-route. Equations (6) and (7) establish that the vehicle is empty when it departs from the depot and the landfill respectively. Equation (8) forces the balance of waste throughout each vehicle route. Equation (9) guarantees that the capacity of the vehicle is not surpassed. Equation (10) indicates that a vehicle only leaves the depot if it has to visit any micro-route. Equation (11) forces the time of the route (micro-route service and traveling time) to not exceed the time limit.

**2.3.2 Mathematical formulation for the Transfer station case**.

The case with transfer station is represented by another model in which we suppose that the fleet of vehicles is based on the transfer station (i.e., the transfer station is also the depot). For the mathematical formulation of this model we used some of the elements already defined in the previous Section 2.2.1. The only additions are the following variables: the binary variable $y_{ij}$ that takes the value 1 if a vehicle uses the path from the exit of micro-route $i \in M'$ and the entrance of micro-route $j \in M'$ and 0 otherwise; the positive continuous variable $u$ that indicates the load of the vehicle after visiting the micro-route $i \in M$; and the positive continuous variable $t$ for the time feasibility constraints, which is associated with the moment of time after passing through the micro-route $i \in M$.

In this way, the following mixed integer programming model is proposed using the formulation of two indices proposed by Miller et al. (1960) in Equations (12) - (21).

$$min \sum_{i,j \in M', \forall j \neq i} y_{ij} d_{ij} \tag{12}$$

Subject to:

$$\sum_{j \in M, \forall j \neq i} y_{ij} = 1, \forall i \in M \tag{13}$$

$$\sum_{j \in M, \forall j \neq i} y_{ji} = 1, \forall i \in M \tag{14}$$

$$u_i - u_j \leq Q(1 - y_{ij}) - q_j, \forall i, j \in M, i \neq j \tag{15}$$

$$q_i \leq u_i \leq Q, \forall i \in M \tag{16}$$

$$t_i - t_j \leq Tt(1 - y_{ij}) - s_j - h_{ij}, \forall i \in M', j \in M, i \neq j \tag{17}$$

$$t_i + h_{i0} \leq Tt(1 - y_{i0}) + Tt, \forall i \in M \tag{18}$$

$$0 \leq t_i \leq Tt, \forall i \in M' \tag{19}$$

$$x \in B \quad u, t \geq 0$$

The proposed objective consists of minimizing the total distance traveled and is expressed in equation (12). Equations (13) and (14) guarantee that each micro-route is visited only once, having a single successor micro-route and a single predecessor micro-route in the path. Equation (15) prevents the formation of sub-tours. Equation (16) guarantees that the capacity of the vehicles is not exceeded. Equations (17), (18) and (19) are maximum time constraints. Equation (20) establishes the binary nature of the variable $x_i$.

## 3. RESULTS AND DISCUSSION

This section presents the implementation details and execution platform, the scenarios that were used for the computational experimentation, and the main results obtained in the experiments.

### 3.1 Implementation details and execution platform.

The input geographic information was processed and visualized using QGIS[3]. Additionally, within QGIS, OpenStreetMap was used for producing the figures of urban areas (OpenStreetMap Contributors, 2017). The mathematical models of section 2 were coded in GAMS version 35.2.0 (GAMS Development Corporation, 2021) and solved with CPLEX Optimization Studio version 20.1 (International Business Machines Corporation, 2021). The resolution platform was a personal computer with 32GB of RAM memory and an Intel(R) Core(TM) i7-4790 CPU @3.60GHz processor.

### 3.2 Description of real-world instances

In this work we consider three cases. The Current Situation (CS) which represents the current system in Bahía Blanca. This comprehends the depot located in the current location of the collection company (Location A of Fig. 1a) and the landfill (Location L of Fig. 1a) as the place in which vehicles unload the waste. Then, we consider two other cases that include a transfer station and only differ in its location. Case Transfer station

---

1 (TS1) corresponds to the location of the transfer station in the current location of the depot (Location A of Fig. 1a), which is approximately 14 km from the landfill and 2 km from the city. This would be a convenient place to locate the transfer stations since it does not require to invest in the leasing of another place. Additionally, since it is the depot, the present collection routes already start and end -after doing an intermediate stop at the landfill to unload waste- in this place. Case Transfer station 2 (TS2) corresponds to the location of the transfer station in an available slot in the industrial park of the city, which is approximately 12 km from the landfill and 7.5 km from the downtown of the city. The industrial park is a special zone in the outskirts of Bahía Blanca planned by the building code of the city for industrial, logistics and related service activities (Location B of Fig. 1a). Companies in this park can benefit from specific tax policies and economies of scale in the provision of services derived from the high concentration of companies. Thus, it was considered a suitable place for locating the transfer station.

| Micro-route Data | | | | | | | |
|---|---|---|---|---|---|---|---|
| Night Shift | | | | Day Shift | | | |
| Micro-route | Micro-route area [km$^2$] | Micro-route distance [km] | Average MSW collection [kg] | Micro-route | Micro-route area [km$^2$] | Micro-route distance [km] | Average MSW collection [kg] |
| 1 | 2.13 | 44.80 | 11,352.00 | 16 | 3.74 | 52.68 | 8,186.00 |
| 2 | 2.56 | 47.29 | 9,285.00 | 17 | 5.81 | 63.74 | 8,997.00 |
| 3 | 1.81 | 42.96 | 10,952.00 | 18 | 3.00 | 53.22 | 10,674.00 |
| 4 | 2.10 | 44.62 | 11,246.00 | 19 | 4.40 | 53.93 | 10,502.00 |
| 5 | 1.92 | 43.61 | 11,183.00 | 20 | 5.29 | 53.87 | 9,796.00 |
| 6 | 2.30 | 45.78 | 10,096.00 | 21 | 10.34 | 91.89 | 7,390.00 |
| 7 | 2.30 | 45.77 | 10,095.00 | 22 | 3.20 | 38.11 | 9,524.00 |
| 8 | 2.52 | 47.09 | 11,811.00 | 23 | 11.51 | 114.81 | 6,601.00 |
| 9 | 2.27 | 45.63 | 10,993.00 | 24 | 11.69 | 102.49 | 10,239.00 |
| 10 | 2.45 | 46.68 | 8,950.00 | 25 | 6.81 | 89.67 | 7,139.00 |
| 11 | 2.17 | 45.02 | 9,254.00 | 26 | 5.19 | 85.03 | 9,565.00 |
| 12 | 2.19 | 45.18 | 11,183.00 | 27 | 8.37 | 74.70 | 7,982.00 |
| 13 | 2.95 | 49.59 | 9,920.00 | 28 | 4.74 | 58.71 | 8,618.00 |
| 14 | 1.83 | 43.04 | 10,245.00 | 29 | 10.85 | 73.38 | 8,240.00 |

| 15 | 2.28 | 46.11 | 12,063.00 | 30 | 3.54 | 55.85 | 11,175.00 |
|----|------|-------|-----------|----|------|-------|-----------|
| - | - | - | - | 31 | 4.99 | 48.20 | 8,733.00 |
| - | - | - | - | 32 | 10.31 | 102.77 | 5,471.00 |

**Table 1**. Area, distance and average waste generation of each micro-route.

As aforementioned, Bahía Blanca has 17 micro-routes in the day shift and 15 micro-routes in the night shift. Table 1 presents for each micro-route: the covered area; the total distance of the typical internal route that is performed by the collection vehicle to gather the waste; and average amount of waste generated by each micro-route. This information was obtained from the database of the waste collection company. In order to analyze the sensitivity of the analysis to the amount of waste we performed different studies in regard to the estimated amount of generated waste. We consider scenarios with 50%, 60%, 70%, 80%, 90% 100%, 105%, and 110% of the estimated amount of waste. This will also help to analyze the competitiveness of the proposed models when subjected to an exceptional rate of waste generation. The travel distances between micro-routes were calculated using Open Source Routing Machine (Luxen and Vetter, 2011) according to the procedure described in Vázquez (2018).

For the collection time of the micro-route, i.e., service time of the micro-route, an empirical formula was devised since the company records were only for the scenarios with a normal amount of waste (the scenario with the 100% of amount of waste). This formula considers that if less waste is to be collected, the average speed of the vehicle inside of the route would be faster and the opposite occurs if more waste is to be collected. Thus, the data of the different amount of waste divided by the travel distance of the internal path of the micro-route and the average speed was approximated by a linear regression as is shown in Fig. 2.

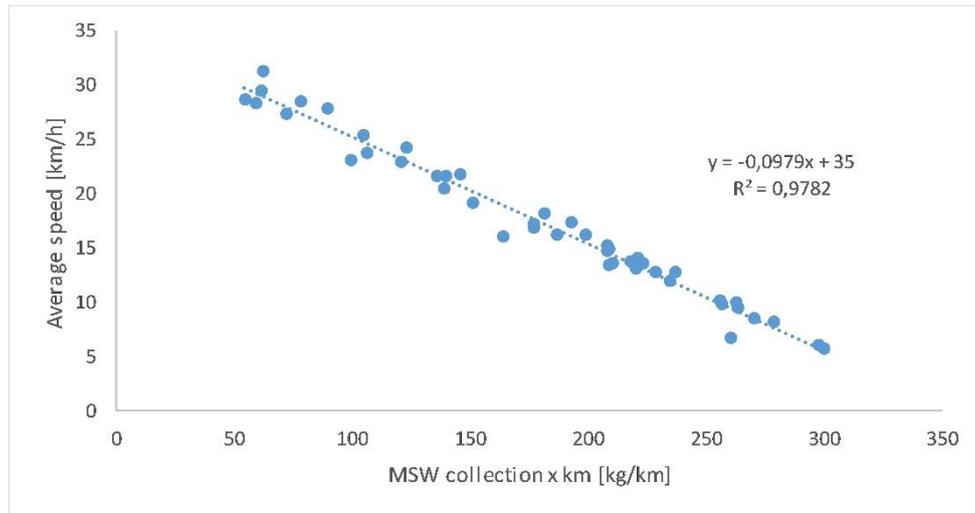

**Figure 2**. Average speed vs. MSW collection x km.

This equation is expressed in Table 2. As an example, in the same table we present the average speed and collection time of the micro-route 22 (depicted in Fig. 1b) in relation to the amount of waste to collect for.

| Micro-route 22 | | | | | | | | |
|---|---|---|---|---|---|---|---|---|
| Data | 50% | 60% | 70% | 80% | 90% | 100% | 105% | 110% |
| Micro-route distance [km] | 38.11 | 38.11 | 38.11 | 38.11 | 38.11 | 38.11 | 38.11 | 38.11 |
| MSW collection [kg] | 4,762.00 | 5,714.40 | 6,666.80 | 7,619.20 | 8,571.60 | 9,524.00 | 10,000.20 | 10,476.40 |
| MSW collection x km [kg/km] | 124.95 | 149.94 | 174.94 | 199.93 | 224.92 | 249.91 | 262.40 | 274.90 |
| Average speed [km/h] | 22.77 | 20.32 | 17.88 | 15.43 | 12.99 | 10.54 | 9.32 | 8.09 |
| Micro-route time [h] | 1.67 | 1.88 | 2.13 | 2.47 | 2.93 | 3.62 | 4.09 | 4.71 |
| Average speed = 35 - 0,0979 x MSW collection x km (empirical equation) | | | | | | | | |

**Table 2**. Average speed and MSW collection time of micro-route 22.

The time limit for the routes was set at eight hours according to the workers shift in the labor convention. The capacity of the trucks of the fleet of vehicles that is currently used in Bahía Blanca is 15,750 kg. These vehicles are used for collecting waste from micro-routes. In the cases in which a transfer station is considered, a larger vehicle with capacity of 25,000 kg is used for taking waste from the transfer station to the landfill. This larger capacity corresponds to commercial vehicles available in the country for compacted waste. As aforementioned, the 32 micro-routes in which the city is partitioned are divided in 17 micro-routes that correspond to the day shift and 15 micro-routes that correspond to the night shift. These shifts are respected in this model and, thus, each scenario is solved in two separated executions: one execution that corresponds in the design of a routing plan for the 17 micro-routes that correspond to the day shift and one execution that corresponds in the design of a routing plan for the 15 micro-routes that correspond to the night shift. The set of input data of the instances can be consulted in the repository of Github[4].

### 3.3 Description of results

In this section we present the main results of the computational experimentation.

**Computing time and CPLEX gap**. Table 3 shows the average computing time and the average gap estimated by CPLEX for both models considering all the runs performed with this model. The gap estimated by CPLEX is indicating the percentage distance between the obtained solution and the estimated optimal solution of the problem (International Business Machines Corporation, 2021). If the solver is able to find the optimal solution, this value would be 0%. However, in computationally complex problems, such as the routing problems that are addressed in this paper (Toth and Vigo, 2002), it is difficult to converge to an optimal result in practical computing times (Rossit et al., 2021). Significant differences in the execution of each model are observed in the values. The CS cases are much harder to solve, requiring, on average, a larger computing time and obtaining solutions that are further from the estimated optimum by

---

[4] https://github.com/diegorossit/Set-of-instances---Molfese-et-al.-2022---WM-R.git

CPLEX. This is in line with the larger complexity of the mathematical model of the CS case (Eqs. (1)-(11)) which has a larger number of restrictions and variables than the TS cases (Eqs. (12)-(19)) due to the inclusion of the landfill as the intermediate facility and, thus, the usage of the three-index formulation of VRP.

| Instance | Average PC runtime [h] | | Average CPLEX gap [%] | |
|---|---|---|---|---|
| | Day shift | Night shift | Day shift | Night shift |
| CS | 1.753 | 0.689 | 24.00 | 26.00 |
| TS1 and TS2 | 0,006 | 0,013 | 0.00 | 0.00 |

**Table 3**. Average PC runtimes y el average CPLEX gap for both models.

**Overall analysis of traveled distance and number of vehicles**. Table 4 summarizes the main information of all the studies performed. From left to right, for each percentage of waste we report, the amount of waste collected, the total traveled distance and number of collection vehicles used in the three cases. Although the number of vehicles is not part of the optimized function in the models of Section 2.2, it is reported since it is a parameter that is generally considered in routing problems (Rossit et al., 2019). To report the total traveled distance, the values of the two shifts are summed and to report the number of collection vehicles, the maximum of the values (of the both shifts) is used. In this Table the vehicles for the trips between the transfer stations and the landfill are not considered. Moreover, the values of traveled distance and the number of vehicles of TS1 and TS2 are reported in percentage difference from the CS case according to Equation (20):

$$(Value_{CS} - Value_{TS1 \, o \, TS2}) * 100/Value_{CS} \qquad (20)$$

where $Value_{CS}$ is the value of the metric considered ( total traveled distance or number of collection vehicles) in the CS case and $Value_{TS1/TS2}$ is the value of the same metric in the transfer station case that is evaluated (TS1 or TS2). Thus, a positive percentage means that there is a reduction of the metric in comparison to the CS.

| Scenario of amount of waste | Total waste collected [kg] | Total traveled distance | | | Number of collection vehicles | | |
|---|---|---|---|---|---|---|---|
| | | CS [km] | TS1 [% diff from CS] | TS2 [% diff from CS] | CS [no] | TS1 [% diff from CS] | TS2 [% diff from CS] |
| 50% | 153,730.00 | 2,268.91 | 2.65% | -1.67% | 7 | 4.29% | 14.29% |
| 60% | 184,476.00 | 2,342.08 | 3.69% | -1.20% | 7 | 0.00% | 0.00% |
| 70% | 215,222.00 | 2,382.36 | 3.80% | -1.22% | 9 | 11.11% | 11.11% |
| 80% | 245,968.00 | 2,495.39 | 5.21% | -1.66% | 9 | 11.11% | 0.00% |
| 90% | 276,714.00 | 2,583.14 | 5.61% | -2.11% | 9 | 0.00% | 0.00% |
| 100% | 307,460.00 | 2,673.45 | 6.05% | -1.34% | 12 | 25.00% | 8.33% |
| 105% | 322,833.00 | 2,711.48 | 7.19% | -0.44% | 14 | 21.43% | 7.14% |
| 110% | 338,206.00 | 2,747.59 | 7.63% | -0.07% | 15 | 20.00% | -6.67% |
| Average % diff from CS | - | - | 5.23% | -1.21% | - | 12.87% | 5.94% |

**Table 4**. Total traveled distance and number of collection vehicles for each percentage of amount of MSW.

The routing plan of the case TS1, which includes the transfer station in an urban location (Location A), has a smaller total traveled distance compared to the CS case in all the scenarios (with an average improvement of 5.23%). However, locating the transfer plant on the outskirts of the city (TS2) has the opposite effect of increasing the total travel distance (with distance on average -1.21% larger than CS case). Both transfer stations cases obtain better results in terms of comparison to CS when the amount of waste is larger. The best percentage improvement from the CS is obtained by TS1 for scenario 110% (7.63%) while the worst is obtained by TS2 in scenario 50% (-1.67%). This result is also explained graphically in Fig. 3a. In the case of the number of the collection vehicles, in general the cases with transfer stations require a smaller number of vehicles except for scenarios with 60%, with 90% and with 80% (only for case TS2). When comparing the number of vehicles between the two alternatives, TS2 uses the same (scenarios 50%, 60%, 70% and 90%) or more vehicles (in the rest of the scenarios) than

TS1. To illustrate these results in Fig. 3 the total traveled distance (a) and number of vehicles (b) are depicted for all the scenarios. For the sake of clarity, the axis of Fig. 3a starts at 2000 km. Although that as far as we are concerned this is the first study of this type in Argentina, the reduction of traveled distances obtained for the TS2 case and number of collection vehicles for the TS1 and TS2 cases in comparison to the CS case is in line with similar studies performed in other countries (Höke and Yalcinkaya, 2021; Rathore and Sarmah, 2019; Vargas et al., 2022).

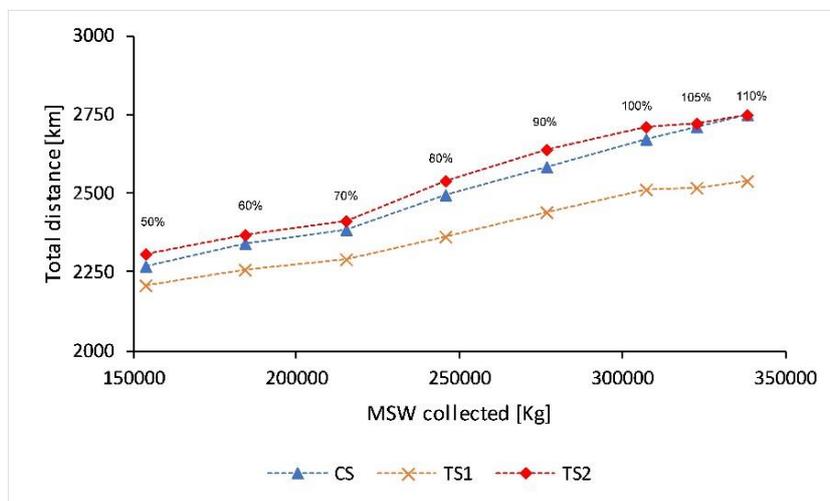

**a)** Total traveled distance.

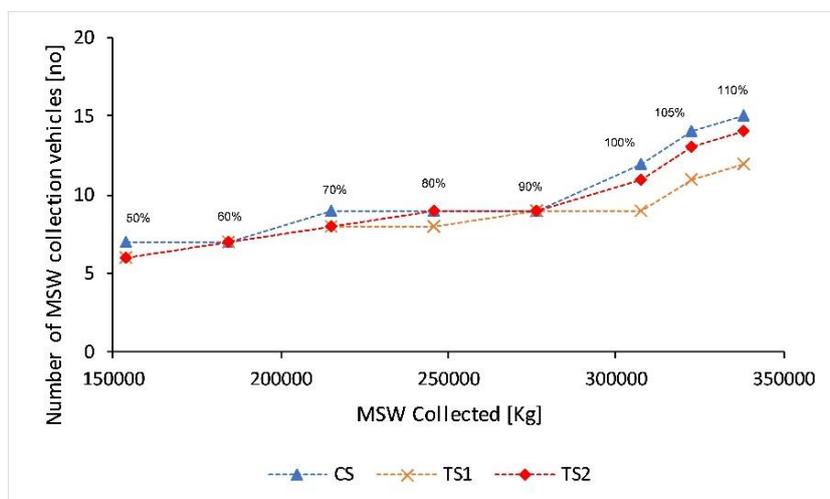

**b)** Number of collection vehicles.

**Figure 3**. Data for each analyzed scenario of amount of waste.

**Specific analysis of the normal scenario**. In Table 5, we present the data from each route that is performed in each routing plan for the three cases and the normal scenario of 100% of amount of waste. This value of the amount of waste is based on the average reported for each micro-route by the collection company of the city. We consider a "route" as the path a collection vehicle performs starting from the depot, visiting a sequence of micro-routes to collect the waste (and the landfill in the case of the CS) and returning to the depot. In Table 5 we also use letters A, B and L for referring to the locations depicted in Fig. 1a. These are Location A (which corresponds to the depot in CS case and the depot and transfer station in TS1), Location B (which corresponds to the depot and transfer station in TS2) and Location L of the landfill. For each route, we present the sequence of stops (micro-routes, transfer stations and/or landfill) that are included, the total distance of the route, the total time of the route and the amount of waste collected. For example, in route 1 of the day shift in the CS case, the sequence of stops is A-8-L-A which means that the vehicle starts at Location A (depot) and then visits the micro-route 8, the Location L (landfill to unload waste) and ends the route at Location A (depot). This route 1 has a total distance of 72.98 km, the collection time is 5.25 h, and the amount of MSW to collect is 11,811 kg.

| Instance | Shift | Route | Sequence of stops | Route distance [km] | Route time [h] | MSW collection [kg] |
|---|---|---|---|---|---|---|
| CS 100% | Night | 1 | A-8-L-A | 72.98 | 5.25 | 11,811.00 |
| | | 2 | A-12-L-A | 79.01 | 5.16 | 11,183.00 |
| | | 3 | A-1-L-A | 71.67 | 5.16 | 11,352.00 |
| | | 4 | A-14-L-A | 74.97 | 4.59 | 10,245.00 |
| | | 5 | A-11-L-6-L-A | 142.99 | 7.93 | 19,350.00 |
| | | 6 | A-9-L-A | 81.55 | 5.02 | 10,993.00 |

| | | 7 | A-10-L-7-L-A | 142.41 | 7.72 | 19,045.00 |
| | | 8 | A-15-L-A | 82.76 | 5.96 | 12,063.00 |
| | | 9 | A-4-L-A | 73.92 | 5.16 | 11,246.00 |
| | | 10 | A-13-L-2-L-A | 147.62 | 7.66 | 19,205.00 |
| | | 11 | A-5-L-A | 72.46 | 5.23 | 11,183.00 |
| | | 12 | A-3-L-A | 73.70 | 5.15 | 10,952.00 |
| | Day | 1 | A-29-21-L-A | 195.48 | 7.31 | 15,630.00 |
| | | 2 | A-26-L-20-L-A | 184.57 | 7.98 | 19,361.00 |
| | | 3 | A-24-L-A | 124.04 | 4.68 | 10,239.00 |
| | | 4 | A-16-L-30-L-A | 166.49 | 7.94 | 19,361.00 |
| | | 5 | A-31-23-L-A | 193.54 | 7.57 | 15,334.00 |
| | | 6 | A-28-32-L-A | 199.35 | 7.38 | 14,089.00 |
| | | 7 | A-18-L-19-L-A | 145.93 | 7.95 | 21,176.00 |
| | | 8 | A-25-27-L-A | 202.16 | 7.42 | 15,121.00 |
| | | 9 | A-22-L-17-L-A | 145.86 | 7.88 | 18,521.00 |
| TS1 100% | Night | 1 | A-1-A-13-A | 106.03 | 7.94 | 21,272.00 |
| | | 2 | A-2-A-5-A | 102.79 | 7.74 | 20,468.00 |
| | | 3 | A-3-A-11-A | 108.24 | 7.88 | 20,206.00 |
| | | 4 | A-4-A-10-A | 110.32 | 7.74 | 20,196.00 |
| | | 5 | A-6-A-9-A | 111.95 | 8.00 | 21,089.00 |
| | | 6 | A-7-A-14-A | 109.48 | 7.68 | 20,340.00 |
| | | 7 | A-8-A | 53.20 | 4.68 | 11,811.00 |
| | | 8 | A-12-A | 57.75 | 4.55 | 11,183.00 |
| | | 9 | A-15-A | 57.07 | 5.31 | 12,063.00 |
| | Day | 1 | A-16-A-26-A | 162.89 | 6.93 | 17,751.00 |
| | | 2 | A-17-A-22-A | 125.10 | 7.29 | 18,521.00 |
| | | 3 | A-19-A-28-A | 132.80 | 6.80 | 19,120.00 |
| | | 4 | A-20-A-29-A | 151.20 | 6.87 | 18,036.00 |
| | | 5 | A-18-A-30-A | 129.54 | 7.67 | 21,849.00 |

| | | | | | |
|---|---|---|---|---|---|
| | | 6 | A-23-21-A | 223.86 | 7.79 | 13,991.00 |
| | | 7 | A-25-27-A | 181.29 | 6.82 | 15,121.00 |
| | | 8 | A-32-31-A | 168.24 | 6.73 | 14,204.00 |
| | | 9 | A-24-A | 114.38 | 4.40 | 10,239.00 |
| TS2 100% | Night | 1 | B-2-B-6-B | 123.22 | 7.27 | 19.381.00 |
| | | 2 | B-5-B-11-B | 108.46 | 8.00 | 20,437.00 |
| | | 3 | B-10-B-13-B | 130.26 | 7.06 | 18,870.00 |
| | | 4 | B-4-B | 63.03 | 4.84 | 11,246.00 |
| | | 5 | B-7-B | 70.65 | 4.12 | 10,095.00 |
| | | 6 | B-1-B | 64.07 | 4.94 | 11,352.00 |
| | | 7 | B-3-B | 66.84 | 4.96 | 10,952.00 |
| | | 8 | B-8-B | 66.02 | 5.05 | 11,811.00 |
| | | 9 | B-9-B | 62.42 | 4.48 | 10,993.00 |
| | | 10 | B-12-B | 65.96 | 4.79 | 11,183.00 |
| | | 11 | B-14-B | 67.96 | 4.39 | 10,245.00 |
| | | 12 | B-15-B | 68.85 | 5.65 | 12,063.00 |
| | Day | 1 | B-16-B-24-B | 178.29 | 7.39 | 18,425.00 |
| | | 2 | B-18-B-22-B | 109.19 | 7.59 | 20,198.00 |
| | | 3 | B-17-B-26-B | 176.77 | 7.35 | 18,562.00 |
| | | 4 | B-19-B-30-B | 127.95 | 7.53 | 21,677.00 |
| | | 5 | B-20-B | 65.80 | 3.47 | 9,796.00 |
| | | 6 | B-21-29-B | 193.20 | 7.24 | 15,630.00 |
| | | 7 | B-25-27-B | 191.89 | 7.13 | 15,121.00 |
| | | 8 | B-28-32-B | 180.81 | 6.85 | 14,089.00 |
| | | 9 | B-31-23-B | 185.24 | 7.34 | 15,334.00 |

**Table 5**. Data of each route for each of the instances in the current situation.

Table 6 presents summary data of the three cases with the normal amount of waste generation. For the cases where a transfer station is used, the number of trips and the distance that must be made to carry the accumulated waste from the station to the landfill per shift are included. The total travel distance and the number of collection vehicles are also summarized for each case in the Table. For example, in the night shift for the TS1 case, the collection vehicle gathers 158,628 kg of waste. Considering that the larger vehicle that performs the trip from the station to the landfill has a capacity of 25,000 kg, the number of trips to carry the waste from the transfer station to the landfill is seven. Since the distance of a round trip from the transfer station to the landfill and back to the transfer station is 25.31 km, the total distance covered in these seven trips is 177,177 km. If the distances of the different routes are added to this distance, the total distance for this shift is 993.99 km. If for this case we consider both shifts, the total distance to travel is 2,511.74 km. Additionally, the number of collection vehicles is defined by the maximum number of routes between both shifts to collect the waste of the micro-routes (the larger vehicle that performs the trip between the transfer station and the landfill is not counted). In Table 6, in TS1 case both shifts require 9 vehicles. However, in the TS2 case the day shift requires 12 vehicles and the night shift requires 9 vehicles. Thus, the maximum number of shifts is 12.

| Instance | Shift | Number of trips from transfer station to landfill [no] | Distance traveled from from transfer station to landfill [km] | Total distance traveled per shift [km] | Total distance traveled [km] | Number of collection vehicles [no] |
|---|---|---|---|---|---|---|
| CS 100% | Night | 0.00 | 0.00 | 1,557.41 | 2,673.45 | 12.00 |
|  | Day | 0.00 | 0.00 | 1,116.04 |  |  |
| TS1 100% | Night | 7.00 | 177.17 | 993.99 | 2,511.74 | 9.00 |
|  | Day | 6.00 | 128.46 | 1,517.76 |  |  |
| TS2 100% | Night | 7.00 | 199.01 | 1,156.76 | 2,709.29 | 12.00 |
|  | Day | 6.00 | 143.40 | 1,552.53 |  |  |

**Table 6**. Summary data of the scenarios with normal waste generation rate.

Examples of routes for these three situations are shown in Fig. 4 below. For the first case (Fig. 4a), on the daytime route, the collection vehicle begins its journey at location A, connecting with micro-routes 25, 22, and 21. Then, it heads to location L for the final disposal of the collection vehicle. trash. Finally, it returns to location A. On the night route, it connects with micro-routes 10 and 6. For the case TS1 (Fig. 4b), on the daytime route, the collection vehicle begins the journey at location A, connecting with micro-routes. - routes 25 and 20. Then, it returns to location A. The garbage is moved to location L for its final disposal by another larger capacity vehicle. On the night route, it connects with micro-routes 9 and 14. For the TS2 case (Fig. 4c), on the daytime route, the collection vehicle begins the journey at location B, connecting with micro-routes 23, 21, and 20. Then, it returns to location B. The garbage is moved to location L for final disposal by another vehicle with greater capacity. On the night route it connects with micro-routes 11, 5 and 8.

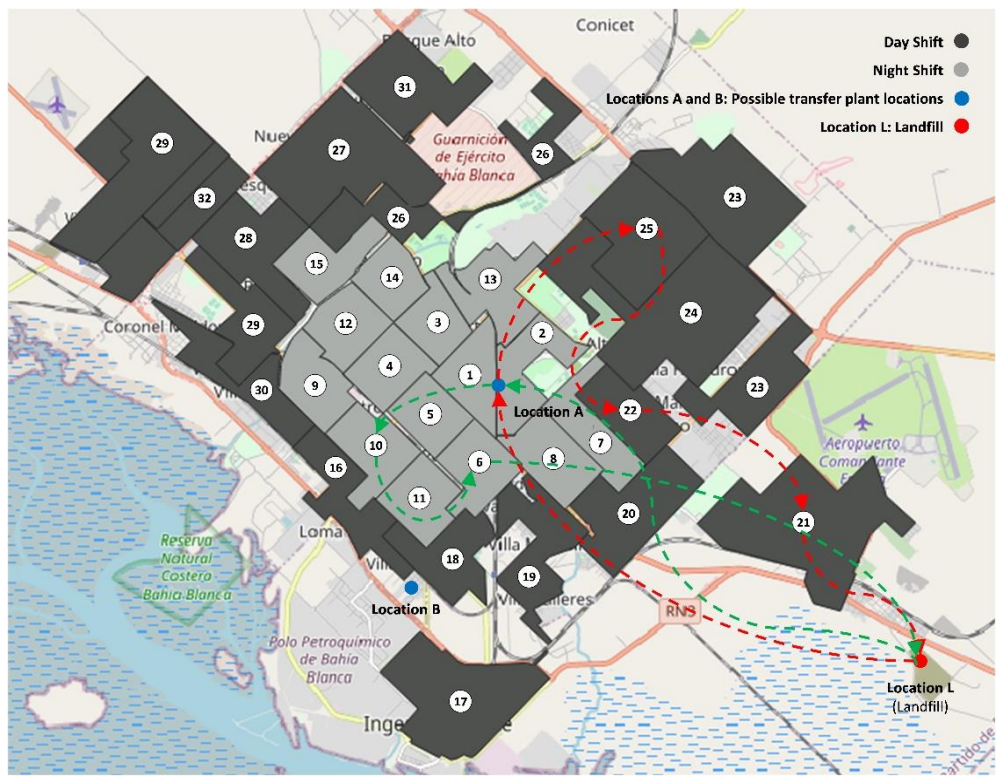

**a)** Routes CS: A-25-22-21-L-A (day shift) / A-10-6-L-A (night shift).

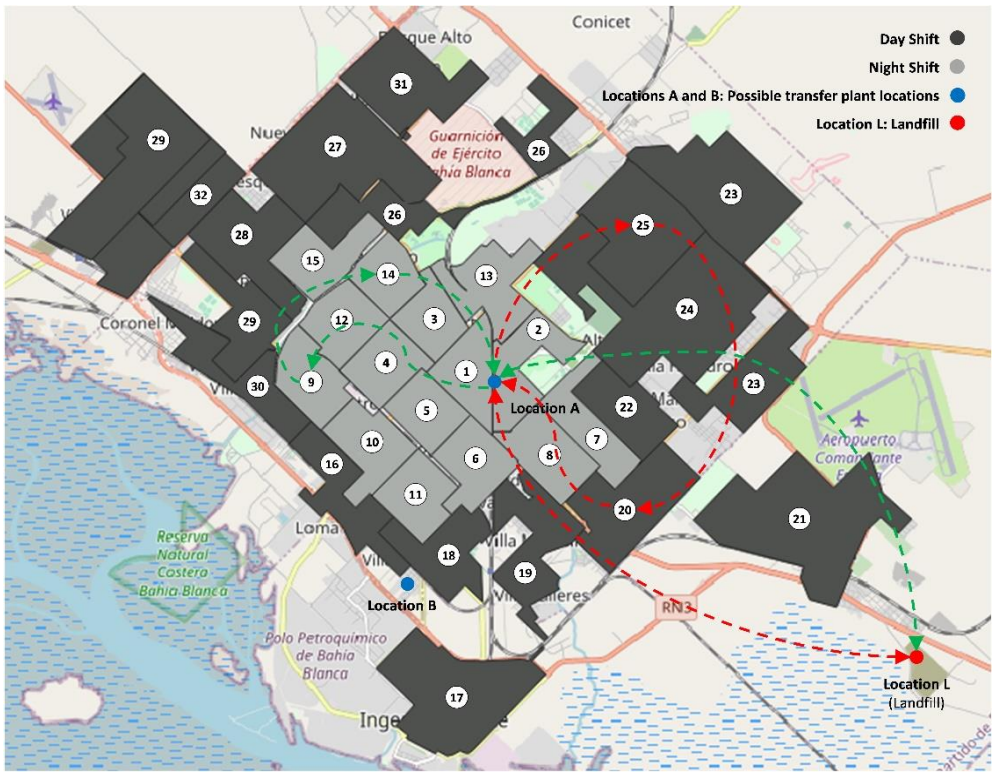

**b)** Routes TS1: A-25-20-A (day shift) / A-9-14-A (night shift).

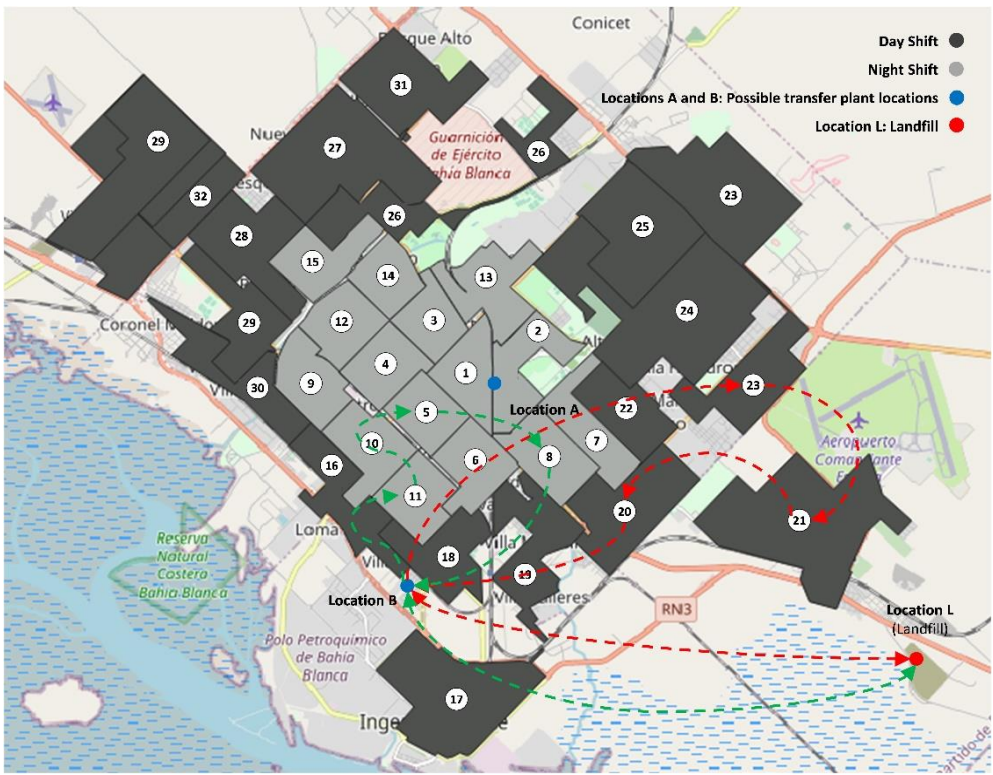

**c)** Routes TS2: B-23-21-20-B (day shift) / B-11-5-8-B (night shift).

**Figure 4.** Examples of routes for CS, TS1 and TS2 (day shift and night shift).

## 4. CONCLUSIONS AND DISCUSSION

Municipal solid waste management is a key issue in modern cities. Not only because when mishandled numerous environmental problems can arise but also because it represents a large proportion of the budgetary expenses of local governments. Among the diverse steps that are involved in this reverse logistic chain, waste collection is one of the most important. In this sense, the present study presents a methodology to optimize waste collection through the sequencing of micro-routes. This methodology considers different operation constraints such as the maximum allowable travel time for a route, that must respect the working shift of the personnel, or the capacity of the collection vehicles that are used in the city.

The proposed methodology allowed to design a waste collection plan for the whole city that constitutes the case study, i.e., Bahía Blanca, Argentina, within reasonable computing times. Moreover, considering that Bahía Blanca do not have a transfer station currently, the proposed methodology is used to analyze the convenience of locating a transfer station. For analyzing this, initially the current situation was studied which comprehends the depot in which collection vehicles are stored and the landfill which is used for unload the waste at the end of each route. This current situation is compared to two scenarios that include transfer station in two different areas of the cities that were selected after a careful qualitative analysis.

The quantitative results show that the transfer station, when located in the current facility of the waste collection company (Location A), is able to reduce the traveled distance on average in 5.25% and up to 7.63% in comparison to the CS case in which no transfer station is used. This is especially important in countries with high logistic costs such as Argentina. The location of the transfer plant in the industrial area of the city (Location B) is not convenient in terms of traveled distance (in comparison to the CS case distances increase on average 1.21% and up to 2.11% in the worst scenario). As qualitative results, the proposed methodology was proved to be adequate for designing the waste collection of a whole city in which micro-routes are already set. This is common in cities in which collection companies already have predefined areas or neighborhoods that due to operational reasons they prefer to do the waste collection all together. Thus, this model can be useful for organizing waste collection in future urban development of the city and even it can be applied to cities that use a similar system based on micro-routes.

The future work main lines are connected to extending the analysis of the micro-routes while considering other aspects in the model, such as the minimization of the number of vehicles, travel times, trips to remote micro-routes and maintenance costs. Additionally, currently the proposed methodology considers that all the accumulated waste is derived to the landfill. This is in line with the current system in Bahía Blanca in which almost all the waste is sent to the landfill. However, in an extended methodology it could be considered that a proportion of the waste is classified –which can be with source classification -at household level- or in the transfer station, and, thus, different recycling centers are included in the model. In this case, the transfer station can work as a classification center and, thus, from this station different trips to other recycling centers can be included. In the case of source classification at household level, the model can be applied directly considering separated tours for each fraction with standard single compartment collection vehicles or modified to include dual compartment collection vehicles. Finally, it would be useful to replicate the sensitivity analysis that was performed in the amount of waste with other parameters of the model, for example the route time limit or the idle capacity of the collection vehicles that is desired by some practitioners to extend the useful life of the collection vehicle.

**ACKNOWLEDGEMENTS**

The first author thanks the Secretariat of Science and Technology of the Universidad Nacional del Sur for the granted scholarship for advanced undergraduate students of the Call 2020.

## REFERENCES

Abarca L, Maas G and Hogland W (2013). Solid waste management challenges for cities in developing countries. *Waste management*, *33*(1), 220-232.

Asefi H, Lim S, Maghrebi M (2015). A mathematical model for the municipal solid waste location-routing problem with intermediate transfer stations. *Australasian Journal* of Information Systems, 19, S21-S35.

Asefi H, Lim S, Maghrebi M and Shahparvari S (2019). Mathematical modelling and heuristic approaches to the location-routing problem of a cost-effective integrated solid waste management. *Annals of Operations Research*, 273(1), 75-110.

Bonomo F, Durán G, Larumbe F and Marenco J (2012). A method for optimizing waste collection using mathematical programming: a Buenos Aires case study. *Waste Management & Research*, *30*(3), 311-324.

Braier G, Durán G, Marenco J and Wesner F (2017). An integer programming approach to a real-world recyclable waste collection problem in Argentina. *Waste Management & Research*, *35*(5), 525-533.

Cavallin A (2019). Análisis de eficiencia y elaboración de propuestas de mejora de la GIRSU en municipios del SO de la Pcia. de Buenos Aires y de Cataluña a través de modelos integrados por DEA y RNA. PhD (Efficiency analysis and proposals for improving the ISWM in municipalities in the SW of the Province of Buenos Aires and Catalonia through integrated models of DEA and ANN). PhD. Thesis. Universidad Nacional del Sur, Argentina.

Cavallin A, Rossit D, Herran V, Rossit, D and Frutos M (2020). Application of a methodology to design a municipal waste pre-collection network in real scenarios. *Waste Management & Research*, *38*(1_suppl), 117-129.

Chatzouridis C and Komilis D (2012). A methodology to optimally site and design municipal solid waste transfer stations using binary programming. *Resources, Conservation and Recycling*, 60, 89-98.

Das S and Bhattacharyya B (2015). Optimization of municipal solid waste collection and transportation routes. *Waste Management*, *43*, 9-18.

Eshet T, Baron M, Shechter M and Ayalon O (2007). Measuring externalities of waste transfer stations in Israel using hedonic pricing. *Waste Management*, *27*(5), 614-625.

Fermani M, Rossit D and Toncovich A (2020). A simulated annealing algorithm for solving a routing problem in the context of municipal solid waste collection. *Communications in Computer and Information Science*, 1408, 63-76.

Gallardo A, Carlos M, Peris M and Colomer F (2015). Methodology to design a municipal solid waste pre-collection system. A case study. *Waste management*, *36*, 1-11.

GAMS Development Corporation (2021). General Algebraic Modeling System (GAMS) Release 35.2.0, Fairfax, USA.

Ghiani G, Manni A, Manni E and Moretto V (2021). Optimizing a waste collection system with solid waste transfer stations. *Computers & Industrial Engineering*, 161, 107618.

Hemmelmayr V, Doerner K, Hartl R and Vigo D (2014). Models and algorithms for the integrated planning of bin allocation and vehicle routing in solid waste management. *Transportation Science*, 48(1), 103-120.

Höke M and Yalcinkaya S (2021). Municipal solid waste transfer station planning through vehicle routing problem-based scenario analysis. *Waste Management & Research*, 39(1), 185-196.

International Business Machines Corporation (2021). *User's Manual for IBM ILOG CPLEX Optimization Studio 22.1.0*. International Business Machines Corporation. Available in: https://www.ibm.com/docs/en/icos/20.1.0 (Accessed: 25-3-2022).

Irnich S, Toth P and Vigo D (2014). Chapter 1: The family of vehicle routing problems. In *Vehicle Routing: Problems, Methods, and Applications*, Second Edition (pp. 1-33). Society for Industrial and Applied Mathematics.

Kirca Ö and Erkip N (1988). Selecting transfer station locations for large solid waste systems. *European Journal of Operational Research*, 35(3), 339-349.

Jia D, Li X and Shen Z (2022). Robust optimization model of waste transfer station location considering existing facility adjustment. *Journal of Cleaner Production*, 340, 130827.

Luxen D and Vetter C (2011). Real-time routing with OpenStreetMap data. In *Proceedings of the 19th ACM SIGSPATIAL* (pp. 513-516).

Markov I, Varone S and Bierlaire M (2015). *The waste collection VRP with intermediate facilities, a heterogeneous fixed fleet and a flexible assignment of origin and destination depot* (TRANSP-OR 150212). Ecole Polytechnique Fédérale de Lausanne, Switzerland. Available in: https://infoscience.epfl.ch/record/208988 (Accessed: 25-3-2022).

Mahéo A, Rossit D and Kilby P (2020). A Benders decomposition approach for an integrated bin allocation and vehicle routing problem in municipal waste management. *Communications in Computer and Information Science*, 1408, 3-18.

Mahéo A, Rossit D and Kilby P (2022). Solving the Integrated Bin Allocation and Collection Routing Problem for Municipal Solid Waste: a Benders Decomposition Approach. *Annals of Operations Research*, in press.

Miller C, Tucker A and Zemlin R (1960). Integer programming formulation of traveling salesman problems. *Journal of the ACM*, 7, 326-329.

Musante C (2021). Indicadores Logísticos Regionales 2021. Technical report. Asociación Latinoamericana de Logística. Available in: https://www.alalog.org/es/studies.

OpenStreetMap Contributors (2017). Planet dump retrieved from https://planet.openstreetmap.org.

Polletto M, Mori P, Schneider V and Zattera A (2016). Urban solid waste management in Caxias do Sul/Brazil: practices and challenges. *Journal of Urban and Environmental Engineering*, 10 (1), 50-56.

Rathore P, and Sarmah S (2019). Modeling transfer station locations considering source separation of solid waste in urban centers: A case study of Bilaspur city, India. *Journal of Cleaner Production*, *211*, 44-60.

Rossit D, Tohmé F, Frutos M and Broz D (2017). An application of the augmented ε-constraint method to design a municipal sorted waste collection system. *Decision Science Letters*, *6*(4), 323-336.

Rossit D, Toutouh J and Nesmachnow S (2020). Exact and heuristic approaches for multi-objective garbage accumulation points location in real scenarios. *Waste Management*, 105, 467-481.

Presentan el presupuesto 2022 en el Concejo Deliberante (2021, December 17). *El Digital de Bahía*. Local.

Rossit D, Vigo D, Tohmé F and Frutos M (2019). Visual attractiveness in routing problems: A review. *Computers & Operations Research*, 103, 13-34.

Rossit D, Nesmachnow S and Toutouh J (2018). Municipal solid waste management in smart cities: facility location of community bins. *Communications in Computer and Information Science*, 978, 102-115.

Rossit D, Toncovich A and Fermani M (2021). Routing in waste collection: A simulated annealing algorithm for an Argentinean case study. *Mathematical biosciences and engineering: MBE*, *18*(6), 9579-9605.

Rossit D and Nesmachnow S (2022). Waste bins location problem: A review of recent advances in the storage stage of the Municipal Solid Waste reverse logistic chain. *Journal of Cleaner Production*, 130793.

Toth P and Vigo D (2002). *The vehicle routing problem*. Society for Industrial and Applied Mathematics.

Toutouh J, Rossit D and Nesmachnow S (2018). Computational intelligence for locating garbage accumulation points in urban scenarios. *Lecture Notes in Computer Science*, 11353, 411-426.

Toutouh J, Rossit D and Nesmachnow S (2020). Soft computing methods for multiobjective location of garbage accumulation points in smart cities. *Annals of Mathematics and Artificial Intelligence*, *88*(1), 105-131.

Vargas A, Díaz D, Jaramillo S, Rangel F, Villa, D and Villegas, J. (2022). Improving the tactical planning of solid waste collection with prescriptive analytics: a case study. *Production*, *32*.

Vázquez Brust A (2018). Ruteo de alta performance con OSRM. *Rpubs by RStudio*. Available in: https://rpubs.com/HAVB/osrm (Accessed: 25-3-2022).

Vazquez Y, Barragán F, Castillo L and Barbosa, S. (2020). Analysis of the relationship between the amount and type of MSW and population socioeconomic level: Bahía Blanca case study, Argentina. *Heliyon*, 6(6), e04343.

Yadav V and Karmakar S (2020). Sustainable collection and transportation of municipal solid waste in urban centers. *Sustainable Cities and Society*, 53, 101937.